\documentclass{article}
\usepackage[english]{babel}
\usepackage{amsfonts, amsmath, amsthm, amssymb,amscd,indentfirst}
\usepackage{amsmath,amssymb,latexsym,indentfirst}
\usepackage{times}
\usepackage{palatino}

\newtheorem{theorem}{Theorem}[section]

\newtheorem{proposition}{Proposition}[section]

\newtheorem{corollary}{Corollary}[section]

\newtheorem{remark}{Remark}[section]

\def\Ac{{\cal A}}

\def\Ac{{\cal A}}

\begin{document}

\title{A rigidity result for the graph case of the Penrose inequality}

\author{Levi Lopes de Lima\thanks{Federal University of Cear\'a,
Department of Mathematics, Campus do Pici, Av. Humberto Monte, s/n, 60455-760,
Fortaleza/CE, Brazil. Partially supported by CNPq/Brazil and FUNCAP/CE. {E-mail: {\tt levi@mat.ufc.br.}}}
\and Frederico Gir\~ao\thanks{Federal University of Cear\'a,
Department of Mathematics, Campus do Pici, Av. Humberto Monte, s/n, 60455-760,
Fortaleza/CE, Brazil. {E-mail: {\tt fred@mat.ufc.br.}}}}

\maketitle
\begin{abstract}
In this note we prove a global rigidity result for asymptotically flat, scalar flat Euclidean hypersurfaces with a minimal horizon lying in a hyperplane, under a natural ellipticity condition.
As a consequence we obtain, in the context of the Riemannian Penrose conjecture, a local rigidity result for the family of exterior Schwarzschild solutions (viewed as graphs in Euclidean space).
\end{abstract}

\section{Introduction and statements of the results}
\label{int}

The purpose of this note is to present a global rigidity result for asymptotically flat, scalar flat Euclidean hypersurfaces $M\subset\mathbb R^{n+1}$ with a minimal horizon lying in a hyperplane $P\subset\mathbb R^{n+1}$, under a natural ellipticity condition. The proof uses a uniqueness result, due to Hounie and Leite \cite{HL}, applied to the two-ended, scalar-flat hypersurface $M'$ obtained by reflecting
$M$  across  $P$, and relies on a regularity argument to make sure that $M'$ is of class $C^2$. As a consequence we obtain, in the context of the Riemannian Penrose conjecture, a local rigidity result for the family of exterior Schwarzschild solutions (viewed as  graphs in Euclidean space). Thus, we start our presentation by recalling the present status of this famous conjecture in General Relativity; see \cite{BC} and \cite{M} for recent surveys on this subject.

Let $(M,g)$ be an $n$-dimensional asymptotically flat Riemannian manifold carrying a (possibly disconneted)
compact inner boundary $\Gamma$ which we assume to be  outermost minimal (we then say that $\Gamma$ is a {\em horizon}). If we assume further that the scalar curvature $R_g$ of $g$ is nonnegative, then the conjectured higher dimensional generalization of the famous (Riemannian) Penrose inequality in General Relativity states that
\begin{equation}\label{penrose}
\mathfrak m_g\geq \frac{1}{2}\left(\frac{|\Gamma|}{\omega_{n-1}}\right)^{\frac{n-2}{n-1}},
\end{equation}
where $\mathfrak m_g$ is the ADM mass of $(M,g)$ and $|\Gamma|$ is the $(n-1)$-area of $\Gamma$,
with the equality occurring if and only if $(M,g)$ is a Schwarzs\-child solution.

For $n=3$ this conjecture has been confirmed in the connected case by  Huisken and Ilmanen \cite{HI}
and in general by Bray \cite{B}. If $n\leq 7$, Bray and Lee \cite{BL} proved the conjecture with the assumption that $M$ is spin for the rigidity statement. Even though many partial results have been obtained, the validity of (\ref{penrose}) remains wide open in higher dimensions except for the case of asymptotically flat Euclidean graphs recently investigated by Lam \cite{L} via a nice integration by parts method which subsequently was extended by  the authors \cite{dLG1} to cover a large class of asymptotically flat hypersurfaces in certain  Riemannian manifolds; see also \cite{dLG2}.
However, we remark that this method, which furnishes integral formulae for the mass, does not seem to be well suited to address the rigidity statement in (\ref{penrose}). As a consequence of our main result (Theorem \ref{main} below) we are able to complement Lam's analysis by providing a local rigidity result for the graph representation of the Schwarzschild solution; see Corollary \ref{pop} below.

More precisely, assume that $\mathbb R^{n+1}$ is endowed with rectangular coordinates $(x,x_{n+1})$, $x\in\mathbb R^n$, and let $M\subset \mathbb R^{n+1}$ be an isometrically immersed hypersurface for which there exists a compact subset $K\subset M$ with the property that $M-K$ can be written as a vertical graph associated with a smooth function $f:\mathbb R^n-K_0\to\mathbb R$, $K_0\subset \mathbb R^n$ compact, satisfying
\begin{equation}\label{asymflat}
f_i(x)=O\left(|x|^{-\frac{n}{2}+1}\right),\quad |x|f_{ij}(x)+|x|^2 f_{ijk}(x)=O\left(|x|^{-\frac{n}{2}+1}\right),
\end{equation}
as $|x|\to+\infty$. Here, $f_i=\partial f/\partial x_i$, etc.
Under these conditions, the {\em ADM mass} of $(M,g)$ is defined by
\begin{equation}\label{mass}
\mathfrak m_g=c_n\lim_{r\to+\infty}\sum_{ij}\int_{S_r}\left(g_{ij,j}-g_{jj,i}\right)\nu_idS_r,
\quad c_n=\frac{1}{(n-1)\omega_{n-1}},
\end{equation}
where $g_{ij}=\delta_{ij}+f_if_j$ are the coefficients of $g$ in non-parametric coordinates  and $\nu$ is the outward unit normal to a large coordinate sphere $S_r$ of radius $r$.

As mentioned above,
an integral formula for $\mathfrak m_g$, in the presence of a horizon $\Gamma\subset M$, has been given in \cite{dLG1}; see (\ref{massform}) below. Instead of rederiving the result, we include here the motivation behind the method, to emphasize how simple it is.
Recall that the Schwarzs\-child solution is given by
\begin{equation}\label{rot}
g_m=\left(1-\frac{2m}{r^{n-2}}\right)^{-1}dr^2+r^2h,
\end{equation}
where $m>0$ is a positive parameter, $r> r_m:=(2m)^{1/n-2}$ is the standard radial coordinate in $\mathbb R^n$ and $h$ is the round metric in $\mathbb S^{n-1}$.
Now, a direct computation using (\ref{mass}) gives $\mathfrak m_g=m$ and since  the horizon $\Gamma$ is the round sphere of radius $r=r_m$, it  follows that
\begin{equation}\label{exac}
\mathfrak m_{g_m}=\frac{1}{2}\left(\frac{|\Gamma|}{\omega_{n-1}}\right)^{\frac{n-1}{n-1}},
\end{equation}
that is, equality holds in (\ref{penrose}), as expected. It turns out, however, that (\ref{exac}) admits a deeper explanation which provides the main motivation for the mass formula (\ref{massform}) below.

To see this, first recall that
the metric (\ref{rot}) can be isometrically embedded in $\mathbb R^{n+1}$ as the graph associated to a radial function $u_m=u_m(r)$, $r=|x|\geq r_m$, satisfying $u_m(r_m)=0$ and
\begin{equation}\label{rot2}
\left(\frac{du_m}{dr}\right)^2=\frac{2m}{r^{n-2}-2m}.
\end{equation}
Such a hypersurface, from now on called a {\em Schwarzschild graph}, clearly satisfies
$S_2=0$, where in general $S_k=S_k(A)$ is the $k^{\rm th}$  elementary symmetric function of the eigenvalues of the shape operator $A$ of a hypersurface (the principal curvatures). In other words, $g_m$ is scalar-flat from an intrinsic viewpoint since $2S_2=R_g$ by the Gauss equation.
On the other hand, for {\em any} hypersurface $M\subset \mathbb R^{n+1}$ endowed with a unit normal $N$, a result by Reilly \cite{R} implies that
\begin{equation}\label{reilly}
{\rm div}_M G(A)X=2S_2\Theta,
\end{equation}
where $G(A)=S_1(A)I-A$ is the first Newton tensor, $X$ is the tangential component of $\partial/\partial x_{n+1}$ and $\Theta=\langle N,\partial/\partial x_{n+1}\rangle$.
In particular, if we think of  the Schwarzschild graph as bounded by the horizon $\Gamma$ and the \lq sphere at infinity\rq\, $S_{\infty}$, which is defined as the limit of large coordinate spheres as $r\to+\infty$, we see from the divergence theorem that
\begin{equation}\label{infty}
\int_\Gamma\langle G(A)X,\nu\rangle d\Gamma=\int_{S_{\infty}}\langle G(A)X,\nu\rangle dS_{\infty},
\end{equation}
where $\nu$ is the outward unit normal and the integral in the right-hand side should be thought of as a limit. It turns out that this formula is just another way of writing (\ref{exac}), since the explicit geometry of Schwarzschild graphs easily gives
$$
c_n\int_\Gamma\langle G(A)X,\nu\rangle d\Gamma=\frac{1}{2}\left(\frac{|\Gamma|}{\omega_{n-1}}\right)^{\frac{n-2}{n-1}},
$$
while another direct computation yields
$$
c_n\int_{S_{\infty}}\langle G(A)X,\nu\rangle dS_{\infty}=m.
$$
Thus, (\ref{exac}) merely reflects the vanishing of the total flux of the divergence free vector field $G(A)X$ over the boundary of the Schwarzschild graph.

Now assume more generally that $M\subset \mathbb R^{n+1}$, an asymptotically flat hypersurface as above, is two-sided in the sense that it carries a globally defined unit normal $N$ which we choose so that $N=\partial/\partial x_{n+1}$ at infinity. Assume further that the inner boundary $\Gamma\subset M$ lies in some {\em horizontal} hyperplane $P$ and that $M$ meets $P$ orthogonally along $\Gamma$. We then say that $\Gamma$ is a {\em horizon} since the orthogonality condition clearly implies that $\Gamma\subset M$ is totally geodesic, hence minimal.
By means of a somewhat more involved computation, again starting from (\ref{reilly}), the following formula for the  mass of $(M,g)$, which generalizes (\ref{exac}), has been proved in \cite{dLG1}:
\begin{equation}\label{massform}
\mathfrak m_g=c_n\int_M\Theta R_gdM+c_n\int_{\Gamma}S_1(\Gamma)d\Gamma,
\end{equation}
where $S_1(\Gamma)$ is the mean curvature of $\Gamma\subset P$ with respect to its inward unit normal. We remark that the graph case of (\ref{massform}) was previously given in \cite{L}, where (\ref{reilly}) is derived by an essentially intrinsic computation.

\begin{remark}\label{more}{\rm
The arguments in \cite{L} and \cite{dLG1} actually establish (\ref{massform}) for a more general type of asymptotics at infinity, but here we restrict ourselves to (\ref{asymflat}) in order to
directly apply the symmetry result in \cite{HL}. In any case, we remark that (\ref{asymflat}) is a rather natural requirement since it means precisely that, at infinity, $M$ approaches a  Schwarzs\-child graph.}
\end{remark}

If we assume,  as in \cite{L}, that $M$ is a graph
(so that $\Theta>0$),
$R_g\geq 0$ and $\Gamma\subset P$ is { convex} then we can apply the well-known Aleksandrov-Fenchel inequality, namely,
\begin{equation}\label{af}
c_n\int_{\Gamma}S_1(\Gamma)d\Gamma\geq \frac{1}{2}\left(\frac{|\Gamma|}{\omega_{n-1}}\right)^{\frac{n-2}{n-1}},
\end{equation}
to conclude  that (\ref{penrose}) holds for $M$; see \cite{S} for an account of (\ref{af}). Moreover, if  equality holds then $R_g=0$
and $\Gamma$ is a union of round spheres. But notice that, regarding the rigidity issue in (\ref{penrose}), this is the best one can achieve by means of (\ref{massform}). Our aim here is precisely to provide a further argument leading to a rigidity result for hypersurfaces (not necessarily graphs) under a natural ellipticity condition (Theorem \ref{main} below). In particular, this yields a local rigidity result for the graph representation of the family of Schwarzschild
solutions in the context of the Penrose conjecture (Corollary \ref{pop}  below).

To explain this, consider an asymptotically flat, scalar flat hypersurface $M\subset\mathbb R^{n+1}$ carrying   a horizon $\Gamma\subset P$
and
require further that the  {\em ellipticity condition}
\begin{equation}\label{elliptic}
S_3(A)\neq 0
\end{equation}
is satisfied everwhere along $M$.
It is not hard to check that this is equivalent to
${\rm rank}\,A\geq 2$. That this is an ellipticity condition can be seen as follows. Scalar-flat hypersurfaces in $\mathbb R^{n+1}$ are critical points, under compactly supported variations, of a natural variational problem, namely, that associated to the functional $\int_MS_1dM$. In this variational setting, the corresponding Jacobi operator is given by
$$
J={\rm div}_M(G(A)\nabla)-3S_3(A).
$$
It turns out that (\ref{elliptic}) is also equivalent to $G(A)$ being positive or negative definite, which means that  $J$ is  elliptic as a differential operator. The reader is referred to \cite{HL} for the proofs of these facts.

\begin{remark}\label{naturell}
{\rm Notice that the ellipticity assumption (\ref{elliptic}) is a natural one, since
it is straightforward to verify that Schwarzschild graphs satisfy ${\rm rank}(A)=n$, hence meeting the condition.}
\end{remark}

\begin{remark}\label{uniq}
{\rm Up to rigid motions, the two-ended hypersurfaces obtained by reflecting (\ref{rot2}) accross the hyperplane $x_{n+1}=0$ exhaust the class of rotationally invariant scalar-flat hypersurfaces in $\mathbb R^{n+1}$; see \cite{Le}.}
\end{remark}

With this terminology at hand, we can now state our main result.

\begin{theorem}\label{main}
Let $M\subset \mathbb R^{n+1}$ be an asymptotically flat, scalar-flat hypersurface carrying a horizon $\Gamma\subset P$  such that $M\cap P=\Gamma$ and assume further that $S_{3}(A)\neq 0$ everywhere along $M$. Then $M$ is (congruent to) a Schwarzschild graph.
\end{theorem}

\begin{remark}\label{stress}
{\rm We stress that regarding the horizon we only assume that $M$ meets $P$ orthogonally along $\Gamma$. In particular, no further geometric restriction (connectedness, convexity, etc.) on the embedding $\Gamma\subset P$ is required. Moroever, the assumption $M\cap P=\Gamma$ is certainly a natural one since it holds for Schwarzschild graphs.}
\end{remark}

The ellipticity of Schwarzschild graphs (Remark \ref{naturell}) leads to the following corollary to Theorem \ref{main}, which yields a {\em local} rigidity result for Schwarzschild graphs in the context of the Penrose conjecture.

\begin{corollary}\label{pop}
Let $M\subset \mathbb R^{n+1}$ be an asymptotically flat, scalar-flat graph carrying a horizon $\Gamma\subset P$  and assume
further that $M$ is a sufficiently small $C^2$ perturbation of a Schwarzschild graph which remains asymptotically flat. Then $M$ is a Schwarzschild graph.
\end{corollary}

The proof follows from the obvious fact that the assumptions $M\cap P=\Gamma$ and $S_3\neq 0$ are both preserved under small $C^2$ perturbations.

\begin{remark}\label{wu}{\rm
A few days after  the first version of this note was published on the arXiv, there appeared a paper by Huang and Wu \cite{HW}, where it is shown that graphs for which the equality holds in the Penrose conjecture are necessarily elliptic, thus extending Corollary \ref{pop} and completely characterizing the equality case. Their reasoning even deals with a more general asymptotics than (\ref{asymflat}).  We remark, however, that from the viewpoint of global rigidity, the argument there seems to work only for graphs and does not cover the more general setting of Theorem \ref{main}.}
\end{remark}

We now briefly sketch the proof of Theorem \ref{main}, which is detailed in the next section.  We first reflect $M$ across the hyperplane $P$ containing the horizon $\Gamma$ so as to obtain  a two-ended hypersurface $M'$ which is scalar-flat everywhere except along $\Gamma$, where it is
only $C^{1,1}$ in principle. We then use the ellipticity condition to prove a regularity result (Proposition \ref{regprop}) showing that actually $M'$ is $C^2$ (in fact, smooth) along the horizon. With this information at hand, it is immediate from asymptotic flatness (\ref{asymflat}) that if we let $P$ be determined by $x_{n+1}=0$ then the graph representation of each end of $M'$ has, as $|x|\to+\infty$, an asymptotic expansion of the type
\begin{eqnarray}
v(x) & = & a|x|^{1/2}+a_1 +a_2|x|^{-1/2}+ |x|^{-3/2}\langle x,c\rangle + O\left(|x|^{-3/2}\right), \label{exp3} \\
v(x) & = & a\log |x| + a_1+ |x|^{-2}\langle x,c\rangle + O\left(|x|^{-2}\right), \label{exp4}\\
v(x) & = & a|x|^{-\frac{n}{2}+2}+ a_1 +|x|^{-n/2}\langle x,c\rangle + O\left(|x|^{-n/2}\right) \label{exp5},
\end{eqnarray}
where $a\neq 0$, $c\in\mathbb R^{n+1}$ and $n=3$, $n=4$ and $n\geq 5$, respectively.
In the language of \cite{HL} this means that each end of $M'$ is {\em regular at infinity} (see their Definition 2.2) and their main result then implies that
$M'$ is rotationally invariant, which means that $M$ is a Schwarzschild graph by Remark \ref{uniq}, as desired. This step of the argument, which relies on the Tangency Principle developed in \cite{HL}, uses not only that $M'$ is elliptic but also that it is embedded, which follows from the assumption $M\cap P=\Gamma$.

\begin{remark}\label{determ}
{\rm We note in passing that the parameter $a\neq 0$ appearing in the asymptotic expansions above admits a nice interpretation. In fact, if we compare the computation for the mass in \cite{dLG1} and the flux formula in Propositon 2.3 of \cite{HL},  it turns out that $a^2$ is proportional to the mass of the corresponding asymptotically flat end. This shows that for an end which is regular at infinity, the leading term in its asymptotic expansion is completely determined by the intrinsic invariant $\mathfrak m_g$.}
\end{remark}

\vspace{0.3cm}
\noindent
{\bf Acknowledgements.} The authors thank M. L. Leite for enlightening conversations during the preparation of this paper.

\section{A regularity result for elliptic scalar-flat hypersurfaces and the proof of Theorem \ref{main}}

As explained above, the proof of Theorem \ref{main} involves the consideration of the {\em embedded} $C^{1,1}$ hypesurface $M'$ obtained from our asymptotically flat, scalar-flat hypersurface $M$ after reflection across the hyperplane $P$ containing the horizon $\Gamma$. More precisely, Theorem \ref{main} follows immediately from the symmetry result in \cite{HL} if we are able to show that $M'$ is actually of class $C^2$ along $\Gamma$. Since the argument is local, we fix $p\in \Gamma$ and write locally $M'$ around $p$ as the graph of a $C^{1,1}$ function $u$ defined in a small neighborhood $U$ of the origin $0\in T_pM'$. Choose rectangular coordinates $(y_1,\cdots,y_n)$ in $U$ so that the hypersurface $\Gamma_0\subset U$ defined by $y_n=0$ is such that $u|_{\Gamma_0}$ is the graph representation of $\Gamma$. Notice that $\Gamma_0$ determines a decomposition $U=U^+\cup U^-$, where $U^+$ (respectively, $U^-$) is given by $y_n\geq 0$ (respectively, $y_n\leq 0$). Clearly, $U^+\cap U^-=\Gamma_0$. We also set $u^{\pm}=u|_{U^{\pm}}$. Moreover, we agree on the index ranges $1\leq i,j,\cdots \leq n$ and $1\leq \alpha,\beta,\cdots\leq n-1$.

We now observe that the following properties hold:
\begin{itemize}
 \item The partial derivatives $u^{\pm}_i$ are $C^1$ along $\Gamma_0$ with $u_i^+=u_i^-$ there;
 \item The function $u^\pm$ is $C^2$ on $U^\pm$ and $u^+_{\alpha\beta}=u^-_{\alpha\beta}$ along $\Gamma_0$.
\end{itemize}

These properties entail the following facts. First, the second property implies that as we approach $\Gamma_0$ by interior points of $U^{\pm}$, all second order derivatives $u_{ij}^\pm$ exist in the limit and are continuous on $U\pm$. The point here is to check whether these derivatives agree along $\Gamma_0$ for each $(i,j)$, so that $u$ is indeed $C^2$ on $U$, which implies that $M'$ is $C^2$  by the fact that $p$ is arbitrarily chosen.
We already know that $u^+_{\alpha\beta}=u^-_{\alpha\beta}$ and, moreover, by the content of the first property applied to $u_n$, we see that $u^+_{\alpha n}=u^-_{\alpha n}$
along $\Gamma_0$  as well. Thus we are led with the task of checking whether $u^+_{nn}=u^-_{nn}$ along $\Gamma_0$.

We notice that $M'_{\pm}=u(U^\pm)$ both have a well-defined shape operator, say $A^{\pm}$, with the usual properties (symmetry, etc.) holding up to $\Gamma_0$. We note for further reference that, in nonparametric coordinates,
\begin{equation}\label{nonparshape}
A_{ij}^\pm=B_{ij}^\pm+C_{ij}^\pm,
\end{equation}
where
\begin{equation}\label{nonparshape2}
B_{ij}^\pm=\frac{u^\pm_{ij}}{W},\quad C_{ij}^\pm=-\frac{1}{W^3}\sum_ku^\pm_iu^\pm_ku^\pm_{kj},\quad W=\sqrt{1+|\nabla u^{\pm}|^2}.
\end{equation}

As usual, given a symmetric matrix $\Ac$, we denote by $S_r(\Ac)$ the $r^{\rm th}$ elementary symmetric function of the eigenvalues of $\Ac$. In particular, we set $S_r(u^\pm)=S_r(A^{\pm})$, so that the following property follows from the assumptions of Theorem \ref{main} and the way $M'$ was constructed from $M$:

\begin{itemize}
 \item $u^\pm$ is an elliptic solution of $S_2(u^\pm)=0$ in the sense that $S_3(u^\pm)\neq 0$.
\end{itemize}

The following proposition provides the regularity result we are looking for.

\begin{proposition}\label{regprop}
Under the conditions above, $u^+_{nn}=u^{-}_{nn}$ along $\Gamma_0$. In particular, $M'$ is of class $C^2$.
\end{proposition}

We start the proof by observing that in general $S_r(\Ac)$ is the sum of the principal minors of order $r$ of the symmetric matrix $\Ac$, so that
\begin{equation}\label{minor}
S_2(\Ac)=\sum_{i<j}\left(\Ac^\pm_{ii}\Ac^\pm_{jj}-(\Ac^\pm_{ij})^2\right).
\end{equation}
It then follows from
(\ref{nonparshape}) that
\begin{equation}\label{minor2}
S_2(u^\pm)=S_2(B^\pm)+S_2(B^\pm,C^\pm)+S_2(C^\pm),
\end{equation}
where
\begin{equation}\label{minor3}
S_2(B^\pm,C^\pm)=\sum_{i<j}\left(B^\pm_{ii}C^\pm_{jj}+B^\pm_{jj}C^\pm_{ii}-2B^\pm_{ij}C^\pm_{ij}\right).
\end{equation}
We now observe that, due to the fact that $u_i^\pm(0)=0$, the ellipticity condition implies that the matrix
\begin{equation}\label{ij}
\frac{\partial S_2(u^\pm)}{\partial u^\pm_{ij}}(0)=
\frac{\partial S_2(B^\pm)}{\partial u^\pm_{ij}}(0)
\end{equation}
is positive or negative definite (see \cite{HL} for a clarification of this point), and we claim that this leads to
\begin{equation}\label{diff}
\sum_\alpha u^\pm_{\alpha\alpha}(0)\neq 0.
\end{equation}
To see this we note that
\begin{equation}\label{decom}
S_2(B^\pm)=\frac{1}{W^2}\sum_\alpha\left(u_{\alpha\alpha}^\pm u_{nn}^\pm-(u^\pm_{\alpha n})^2\right)+
\sum_{\alpha<\beta}\left(B_{\alpha\alpha}^\pm B_{\beta\beta}^\pm-(B^\pm_{\alpha \beta})^2\right),
\end{equation}
with the second term on the right-hand side {\em not} depending on $u_{nn}$. It then follows that
$$
\frac{\partial S_2(B^\pm)}{\partial u_{nn}}=\frac{1}{W^2}\sum_\alpha u^\pm_{\alpha\alpha},
$$
and we see that the left-hand side in (\ref{diff}) is precisely the $(n,n)$-entry of (\ref{ij}), so  the claim follows.

We now use that $S_2(u^\pm)=0$ in (\ref{minor2}), which together with (\ref{decom}) gives
\begin{equation}\label{decom2}
D^\pm u_{nn}^\pm+E^\pm+F^\pm=0,
\end{equation}
where
$$
D^\pm=\frac{1}{W^2}\sum_\alpha u^\pm_{\alpha\alpha},
$$
$$
E^\pm=-\frac{1}{W^2}\sum_\alpha (u^\pm_{\alpha n})^2+ \sum_{\alpha<\beta}\left(B_{\alpha\alpha}^\pm B_{\beta\beta}^\pm-(B^\pm_{\alpha \beta})^2\right),
$$
and
$$
F^\pm=S_2(B^\pm,C^\pm)+S_2(C^\pm).
$$
It follows from (\ref{diff}) that $D^\pm$ does not vanish and remains bounded in a neighborhood of the origin, so that
$$
u_{nn}^\pm=-\frac{E^\pm}{D^\pm}-\frac{F^\pm}{D^\pm}.
$$
in this neighborhood. Due to the fact that $F^{\pm}$ depends at least quadratically on the first order derivatives, the last term in the right-hand side vanishes as we approach the origin, since all second order derivatives, including $u^\pm_{nn}$, remain bounded there. On the other hand, the first term in the right-hand side only depends on $u^\pm_i$ and $u^{\pm}_{\alpha i}$ and since they coincide at the origin we conclude that $u_{nn}^+(0)=u_{nn}^-(0)$, as desired.
This completes the proof of Proposition \ref{regprop} and hence, by the comments on the paragraph immediately before Remark \ref{determ}, of Theorem \ref{main}.

\end{document}